\documentclass{amsproc}
\usepackage{amsfonts}
\usepackage{amsopn}
\usepackage{epsfig}
\usepackage{amsmath}
\begin{document}

\newtheorem{theorem}{Theorem}[section]
\newtheorem{proposition}[theorem]{Proposition}
\newtheorem{lemma}[theorem]{Lemma}
\newtheorem{definition}[theorem]{Definition}
\newtheorem{remark}{Remark}[section]
\newtheorem{corollary}[theorem]{Corollary}
\newtheorem{conjecture}[theorem]{Conjecture}

\renewcommand{\P}{\mathbb{P}}
\newcommand{\C}{\mathbb{C}}
\renewcommand{\O}{\mathcal{O}}
\newcommand{\G}{\mathbb{G}}
\renewcommand{\a}{\sum_{i=1}^a \alpha_i}
\renewcommand{\b}{\sum_{j=1}^b \beta_j} 
\newcommand{\X}{\mathcal{X}}
\newcommand{\Y}{\mathcal{Y}}
\newcommand{\Z}{\mathbb{Z}}
\newcommand{\M}{\overline{M}_{0,n}}
\newcommand{\F}{\mathbb{F}}

\title{The Arithmetic and the Geometry of Kobayashi Hyperbolicity}
\author{Izzet Coskun}  \address{Harvard University, Mathematics
  Department, Cambridge, MA 02138} \email{coskun@math.harvard.edu} \date{ } 
\subjclass{Primary 32Q45} 
\thanks{During the preparation of this
    article the author was partially supported by a Clay Mathematics
    Institute Liftoff Fellowship.}
\maketitle

Curves of genus greater than one exhibit arithmetic and geometric
properties very different from curves of genus zero and one.  In these
notes we will survey a few ways of extrapolating these properties to
higher dimensional complex manifolds. We will specifically concentrate
on Brody and Kobayashi hyperbolicity, which are generalizations of
complex analytic properties of higher genus curves. \smallskip

The basic metric and geometric properties that distinguish curves of
genus two or more from curves of genus zero and one can be listed as
follows.

\begin{enumerate}
\item The dimension of the pluricanonical series, $h^0(C,mK)$, grows
  linearly with $m$ for curves of genus at least two.
  
\item The canonical/cotangent bundle of a curve of genus at least two
  is ample.
  
\item A curve of genus at least two admits a hyperbolic metric with constant
  negative curvature.
  
\item Curves of genus at least two are uniformized by the unit disc,
  hence they do not admit any non-constant holomorphic maps from $\C$ .
\end{enumerate}
Each of these properties can be generalized to higher dimensional
manifolds, often in many distinct ways. For example, the growth rate of
the dimension of sections of the pluricanonical series leads to the
concept of Kodaira dimension, a useful and well-studied birational
invariant. \smallskip

In these notes we will discuss generalizations of the last two
properties to higher dimensional varieties keeping in mind the Schwarz
Lemma and Liouville's Theorem from the theory of complex
functions. \smallskip

\noindent{\bf Acknowledgements:} I would like to thank L. Caporaso,
S. Grushevsky, J. Harris and R. Vakil for fruitful discussions and the
referee for many useful comments.

\section{Introductory remarks about hyperbolicity}

In this section we define Brody and Kobayashi hyperbolicity and state
their basic properties. We also give examples of hyperbolic and
non-hyperbolic manifolds. The reader can refer to
\cite{lang:hyperbolic} and \cite{demailly:hyperbolic} for more
details.\footnote{After these notes were written, I discovered that O.
  Debarre also has some very nice notes on hyperbolicity
  \cite{debarre:hyperbolicity}. The reader is encouraged to consult
  these notes for additional information on hyperbolic varieties,
  especially about those varieties with ample cotangent bundle.}  Let
$\C$ denote the complex plane. We use $X$ and $Y$ to denote complex
manifolds. We reserve $C$ for curves.  \smallskip

\noindent {\bf Observation:} If $C$ is a curve of genus at least two,
then any holomorphic map $f: \C \rightarrow C$ is necessarily
constant.  \smallskip

\noindent {\bf Proof:}  Any map $f:\C \rightarrow C$ factors through
the universal cover of $C$, which is the unit disc. Such a map is
constant since by Liouville's Theorem every bounded entire function is
constant.  $\Box$ \smallskip

In contrast, curves of genus zero and one do admit non-constant
holomorphic maps from $\C$. This property of higher genus curves can
be generalized to higher dimensional manifolds and leads to the
concept of Brody hyperbolicity.

\begin{definition}
  A complex manifold $X$ is {\bf Brody hyperbolic} if there are no
  non-constant holomorphic maps from $\C$ to $X$.
\end{definition}        
{\bf Examples:} $\bullet$ Any variety of the form $\prod_{i=1}^n C_i$,
where $C_i$ are curves of genus at least two, is Brody
hyperbolic. More generally, a finite product of Brody hyperbolic
manifolds is Brody hyperbolic.
 \smallskip
 
 $\bullet$ If a complex manifold $X$ contains a rational curve or a
 complex torus, then $X$ is not Brody hyperbolic. If $X$ contains
 a rational curve, then the inclusion of $\C$ in the rational curve
 followed by the inclusion of the rational curve in $X$ gives a
 non-constant holomorphic map from $\C$ to $X$.
 
 The universal cover of a $g$-dimensional complex torus is $\C^g$.
 Taking the image of a general complex line in $\C^g$ under the
 quotient map gives a non-constant holomorphic map from $\C$ into the
 complex torus. If $X$ contains a complex torus, composing the
 previous map with the inclusion of the complex torus in $X$ gives a
 non-constant holomorphic map from $\C$ into $X$. \smallskip

 $\bullet$ The blow-up of a variety is not Brody hyperbolic.
 Consequently, Brody hyperbolicity is not a birational invariant.
 \smallskip
 
 $\bullet$ Consider $X = \P^2 - \cup_{i=1}^4 l_i$ where $l_i$ are
 general lines.  Then $X$ is not Brody hyperbolic since there are
 $\C^*$'s in $X$---take the intersection of $X$ with the line passing
 through $l_1 \cap l_2$ and $l_3 \cap l_4$.  \smallskip

\begin{figure}[htbp]
\begin{center}
\epsfig{figure=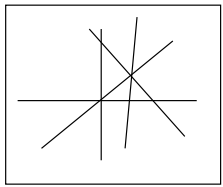}
\end{center}
\caption{A Brody hyperbolic manifold.}
\label{lines}
\end{figure}

$\bullet$ Consider $Y = \P^2 - \cup_{i=1}^5 l_i$ where $l_i$ are the
lines pictured in Figure \ref{lines}.  $Y$ is Brody hyperbolic.
Consider the pencil of lines in $\P^2$ based at one of the points
where three of the lines $l_i$ intersect. Restricting this pencil to
$Y$, we see that $Y$ is fibered over $\P^1$ punctured at three points
with fibers isomorphic to $\P^1$ punctured at three points. Composing
any  holomorphic map from $\C$ with the map to the base of the
fibration, we conclude that the image of the map must lie in a fiber.
Since the fibers are Brody hyperbolic, the map must be constant.
Using the following easy lemma, which generalizes the idea of this
example, one can generate many examples of Brody hyperbolic manifolds.

\begin{lemma}\label{fibration}
  Let $f: X \rightarrow Y$ be a smooth morphism of complex manifolds,
  where $Y$ is Brody hyperbolic and $f^{-1} (y)$ is Brody hyperbolic
  for every $y \in Y$. Then $X$ is Brody hyperbolic.
\end{lemma}

The Schwarz-Pick Lemma  states that a holomorphic
map between two hyperbolic Riemann surfaces is either a local isometry
or distance decreasing. The hyperbolic distance between any two points
$p$ and $q$ on a hyperbolic Riemann surface $C$ is equal to the
shortest hyperbolic distance between lifts of the points $p$ and $q$
to the universal cover, the unit disc $\Delta$. Any holomorphic map
$f: \Delta \rightarrow C$ factors through the universal cover.
Consequently, the Schwarz-Pick Lemma implies that the infimum of the
hyperbolic distances between $p'$ and $q'$ in $\Delta$ for which there
exists a holomorphic map $f: \Delta \rightarrow C$ such that $f(p')=p$
and $f(q') = q$ is achieved when $f$ is the universal covering
map. Moreover, this infimum 
is equal to the hyperbolic distance between $p$ and $q$. In this form
we can generalize the hyperbolic distance to higher dimensional
manifolds. \smallskip

Let $\Delta_r$ denote the disc of radius $r$ in the complex plane.  We
will denote the unit disc by $\Delta$.  If $X$ is a complex manifold,
then $X$ can be endowed with a pseudo-distance due to Kobayashi.

\begin{definition}
  Given a tangent vector $\xi \in T_{X,x}$ at $x \in X$ we define its
  {\it Kobayashi pseudo-norm} to be
  $$k(\xi) := \inf_{\lambda} \{ \lambda : \exists \ f : \Delta
  \rightarrow X, \ f(0) = x, \ \lambda f'(0) = \xi \}.$$
  The {\it
    Kobayashi pseudo-distance} $d_{X}$ is the geodesic pseudo-distance
  obtained by integrating this pseudo-norm.
\end{definition}

\begin{figure}[htbp]
\begin{center}
\epsfig{figure=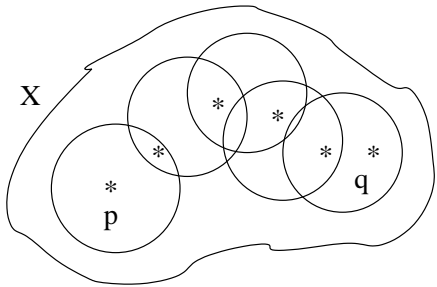}
\end{center}
\caption{Measuring the Kobayashi distance.}
\label{kob}
\end{figure}

More visually, we can describe how to compute the distance between any
two points $p,q \in X$ as follows. We find chains of maps $f_i :
\Delta \rightarrow X$ with two distinguished points $p_i, q_i$, where
$f_1(p_1) = p$, $f_n(q_n)= q$ and $f_i(q_i) = f_{i+1}(p_{i+1})$. The
Kobayashi pseudo-distance is the infimum of the sum of the hyperbolic
distances between $p_i$ and $q_i$ over all such chains. \smallskip

The Schwarz-Pick Lemma generalizes to holomorphic maps between higher
dimensional manifolds endowed with the Kobayashi pseudo-distance.

\begin{lemma}\label{Pick}
  If $f: X \rightarrow Y$ is a holomorphic map between two complex
  manifolds, then $d_X(p,q) \geq d_Y (f(p), f(q))$.
\end{lemma} 
Unfortunately, $d_X$ does not have to be non-degenerate. \smallskip

\noindent {\bf Main Counterexample:} The Kobayashi pseudo-distance on
$\C$ is identically zero. To compute the distance between two points
$x$ and $y$ in $\C$,  consider for each integer $n>1$ the functions
$f_n(z)= n(y-x)z + x$ from 
the unit disc to $\C$. The function $f_n(z)$ maps 0
to $x$ and $1/n$ to $y$. Since the hyperbolic distance between $0$ and
$1/n$ in the unit disc tends to zero as $n$ tends to infinity, we
conclude that the Kobayashi pseudo-distance between $x$ and $y$ in
$\C$ is zero. \smallskip

\begin{definition}
  A complex manifold $X$ is called {\it Kobayashi hyperbolic} if its
  Kobayashi pseudo-distance is non-degenerate.
\end{definition}

\noindent {\bf Kobayashi hyperbolicity implies Brody hyperbolicity. }
  As a consequence of Lemma \ref{Pick} and the above example we see that
if a complex manifold $X$ admits a non-constant holomorphic map from
$\C$, then $X$ cannot be Kobayashi hyperbolic. In other words,
Kobayashi hyperbolicity implies Brody hyperbolicity.  For compact
complex manifolds the converse of this result is true
(\cite{brody:hyperbolic}).

\begin{theorem}[Brody]\label{Br=Ko}
  A compact complex manifold $X$ is Kobayashi hyperbolic if and only
  if it is Brody hyperbolic.
\end{theorem} 

\noindent{\bf Sketch of proof:} If $X$ is not Kobayashi hyperbolic, then there
exists a sequence of maps $f_n : \Delta \rightarrow X$ such that
$|df_n (0 )|$ tends to $\infty$, where $| \cdot |$ denotes a fixed
Hermitian metric on $X$. By rescaling the map we can assume that we
have a sequence of maps $f_n : \Delta_{r_n} \rightarrow X$ where
$|df_n(0)| = 1$ and the radii $r_n$ tend to $\infty$. \smallskip
  
Let $f:\Delta_r \rightarrow X$ be a holomorphic map with $|df(0)| =
c$. Let $f_t(z) = f(tz)$. Consider the function $$g(t)= \sup_{z \in
  \Delta_r} |df_t (z)|.$$
The function $g(t)$ is  increasing
for $0 \leq t \leq 1$ and continuous for $0 \leq t <1$. Moreover, as
$t$ tends to 1 from below $g(t)$ tends to $g(1)$.  Hence there exists
$t$ in the interval $[0,1]$ and an automorphism $h$ of
$\Delta_r$ such that $$\sup_{z \in \Delta_r} |d(f_t \circ h)
(z)| = |d(f_t \circ h) (0)| = c .$$
This is known as the Brody
reparametrization lemma. \smallskip
  
Applying the Brody reparametrization lemma to our family of maps $f_n$, we
 obtain
 a new family of
maps $\tilde{f}_n :\Delta_{r_n} \rightarrow X$ whose derivatives are bounded
in norm by 1 in $\Delta_{r_n}$. Moreover, $|d \tilde{f}_n (0)| = 1$
 for every $n$. We would like to show that we can select a subsequence
 from this family that is uniformly convergent on every compact subset
 of $\C$. Since the derivatives at $0$ all have norm 1, then the
 family converges to a non-constant holomorphic map from $\C$ to
 $X$. This violates Brody hyperbolicity. \smallskip

Using the compactness of $X$ we can check uniform convergence in a
neighborhood of every point $z_0$ in $\C$. Again using the compactness
of $X$ (and passing to a subsequence if necessary) we can assume the
family converges at $z_0$. In this situation the derivative bounds
imply that the  family forms an equicontinuous family of
holomorphic maps around $z_0$.  We can conclude that a subsequence
converges uniformly on compact sets by Ascoli's theorem. We thus
obtain a holomorphic map from $\C$ to $X$. $\Box$

\smallskip

Note that Brody's theorem does not have to hold for non-compact
manifolds (see \cite{brody:thesis}).  Consider the domain in $\C^2$
given by
$$
D : = \{ (z,w) \ | \ |z|<1, \ |zw|<1 \ \mbox{and} \ |w| < 1 \ 
\mbox{if} \ z = 0 \}. $$
The projection of the domain $D$ to the first
coordinate gives rise to a family of discs parameterized by the unit
disc in the $z$-plane.  Consequently, Lemma \ref{fibration} implies
that $D$ is Brody hyperbolic. On the other hand, when $z$ approaches
zero, the radius of the disc lying over $z$ tends to infinity. Using
this we can show that the Kobayashi pseudo-distance between $p =
(0,0)$ and $q =(0,w_0)$ in $D$ vanishes.  We can find two points
$p'$ and $q'$ having the same $w$ coordinates as $p$ and $q$,
respectively, in a fiber close to zero. Since $p'$ and $q'$ are in a
very large disc, the hyperbolic distance between them is very small.
Since the hyperbolic distances between $p$ and $p'$ and $q$ and $q'$
are very small, the sum of the hyperbolic distances can be made
arbitrarily small. \smallskip

More explicitly, consider  the three maps
$$f_1 (z) = (z,0), \ f_2(z) = \left( \frac{1}{n}, nz \right) , \ f_3
(z) = \left( \frac{1}{n} + \frac{1}{2} z, w_0 \right)$$
from the unit
disc into $D$. Note that $f_1(0)= (0,0)$, $f_1(1/n)=(1/n,0)$, $f_2(0)=
(1/n,0)$, $f_2(w_0/n) = (1/n, w_0)$ and $f_3(0)= (1/n,w_0)$,
$f_3(-2/n)= (0,w_0)$.  Hence these maps form a chain of maps from the
unit disc to $D$ connecting $(0,0)$ to $(0,w_0)$. As $n$ tends to
infinity, the sum of the hyperbolic distances between the chosen
points in the unit disc tends to zero. We conclude that the Kobayashi
pseudo-distance between $(0,0)$ and $(0,w_0)$ is zero. Observe that
this example also shows that the analogue of Lemma \ref{fibration}
does not hold for Kobayashi hyperbolicity.  \smallskip

Under suitable hypotheses on a subvariety $Y$ of $X$, we can assert
that $X-Y$, the complement of $Y$ in $X$, is Kobayashi hyperbolic. For
example, the following theorem is often useful in proving that the
complement of a subvariety is Kobayashi hyperbolic.

\begin{theorem}
  Let $X$ be a compact variety and let $Y$ be a proper algebraic subset. If
  $Y$ and $X - Y$ are Brody hyperbolic, then $X-Y$ is Kobayashi
  hyperbolic.
\end{theorem}

\noindent {\bf Sketch of proof:} If $X-Y$ is not Kobayashi hyperbolic,
we can get a sequence of maps $f_n: \Delta_{r_n} \rightarrow X-Y$ that
converges to a holomorphic map $g$ from $\C$ to $X$.  The question is
whether the image of $g$ intersects $Y$.  The image cannot be
contained in $Y$ since $Y$ is Brody hyperbolic.

If we rule out that the image of $g$ intersects $Y$, then we obtain a
contradiction since $X-Y$ is Brody hyperbolic. One can prove that the
intersection points of the image of $g$ with $Y$ have to be isolated.
Suppose $g(z_0) \in Y$. Take a circle $S$ around $z_0$ such that $g(S)
\subset X-Y$. The winding number is zero because $f_n(\Delta_{r_n})$
does not meet $Y$. Hence $Y$ cannot intersect the image of the
interior of $S$, contradicting that $g(z_0) \in Y$. $\Box$ \smallskip

We now discuss some examples of Kobayashi hyperbolic
manifolds. From now on whenever we say hyperbolic without further
qualification, we will always mean Kobayashi hyperbolic.  The
following example due to Green (\cite{green:hyperplane}) gives the
first non-trivial examples.

\begin{theorem}[Green]
  The complement of $2n+1$ general hyperplanes in $\P^n$ is
  hyperbolic.
\end{theorem}

Generalizing further one can ask whether the complement of a very
general irreducible hypersurface in $\P^n$ of large enough degree is
hyperbolic.  Siu and Yeung answer this question positively in $\P^2$
(\cite{siuyeung:complement}). Here and in the following by `very
general' we will refer to the complement of the union of countably
many proper subvarieties.

\begin{theorem}[Siu-Yeung]\label{curve}
  Let $C$ be a very general curve of sufficiently large degree
  in $\P^2$. Then $\P^2 - C$ is Kobayashi hyperbolic.
\end{theorem}

Note that this theorem can be interpreted as a generalization of
Picard's theorem to higher dimensions. Recall that Picard's theorem
says that any entire map which omits two values is constant. This
theorem implies that any holomorphic map from $\C$ into $\P^n$ which
omits a (very general) hypersurface of a large enough degree is constant.
\smallskip

Siu and Yeung also prove the analogue of this result for the
complement of very general ample divisors in abelian varieties
(\cite{siuyeung:4} \cite{siuyeung:one}, \cite{siuyeung:2}
\cite{siuyeung:3}). 

\begin{theorem}[Siu-Yeung]
  Let $A$ be an abelian variety. Let $D$ be a very general ample
  divisor in $A$. Then $A - D$ is Kobayashi hyperbolic.
\end{theorem} 

\noindent  There is a close relation between the hyperbolicity of
complements of divisors in projective space and the hyperbolicity of
general hypersurfaces in projective space. We have the following
theorem due to Siu (\cite{siu:hyperbolic}, \cite{siu:hypersurface}).

\begin{theorem}[Siu]\label{siuhyp}
  A very general surface of sufficiently high degree in $\P^3$ is
  Kobayashi hyperbolic.
\end{theorem}

\noindent In fact, Siu recently has generalized Theorems \ref{curve} and
\ref{siuhyp} to hypersurfaces in $\P^n$. Not all the proofs
have appeared in print.

\begin{theorem}[Siu]
  Let $X$ be a very general hypersurface of sufficiently large degree
  in $\P^n$. Then $\P^n - X$ is Kobayashi hyperbolic.
\end{theorem}  

\begin{theorem}[Siu]\label{higherdim}
  A very general hypersurface of sufficiently high degree in $\P^n$ is
  Kobayashi hyperbolic.
\end{theorem}

\noindent For surfaces in $\P^3$ the term `sufficiently high degree'
can be taken to mean at least 11. Conjecturally a very general
hypersurface of degree larger than $2n$ in $\P^n$ is expected to be
hyperbolic.  However, the current known bounds are much larger than the
expected bounds.  \smallskip

Finally, we observe that some properties of the tangent or cotangent
bundle of a complex manifold forces the manifold to be hyperbolic. For
example, if $X$ admits a metric with negative sectional curvature; or
if $X$ has ample cotangent bundle; or if $X$ is the quotient of a
bounded domain in $\C^n$ by a free group action, then $X$ is Kobayashi
hyperbolic. \smallskip

As an amusing corollary we note that $M_g^0$, the moduli space of
automorphism-free, smooth curves of genus $g > 2$, is Kobayashi
hyperbolic since the Weil-Petersson metric is a metric on $M_g^0$ with
negative sectional curvature.  Consequently, $M_g^0$ does not contain
any complete rational or elliptic curves. It would be interesting to
give lower bounds on the genus of complete curves contained in $M_g^0$.

\section{The geometry of Kobayashi hyperbolicity}

In this section we discuss some proven and conjectural geometric
characterizations of hyperbolicity. Hyperbolicity imposes strong
restrictions on the geometry of a variety. In particular it constrains
the type of subvarieties the variety can have.

\begin{proposition}
  Let $X$ be a compact complex hyperbolic manifold with Hermitian
  metric $\omega$. Then $\exists  \ \epsilon > 0$ such that for any
  reduced, irreducible curve $C \subset X$, the genus $g$ of the
  normalization satisfies $$g \geq \epsilon \deg_{\omega} (C)$$
\end{proposition} 

\noindent {\bf Sketch of proof:} If $X$ is hyperbolic, then there
exists a constant such that $d_X(\xi) \geq \epsilon_0 || \xi
||_{\omega}$ for every tangent vector $\xi$. Let $\nu: C^{\nu}
\rightarrow C$ be the normalization of a curve $C \subset X$. If we
denote the hyperbolic metric on $C^{\nu}$ by $k_{C^{\nu}}$, then the
Gauss-Bonnet formula implies that $$-\frac{1}{4} \int_{C^{\nu}}
\mbox{curv}(k_{C^{\nu}}) = -\frac{\pi}{2} \chi(C^{\nu}).$$
There is a
natural holomorphic map $i \circ \nu: C^{\nu} \rightarrow X$ obtained
by composing the normalization map $\nu$ from $C^{\nu}$ to $C$ by the
inclusion $i$ of $C$ in $X$. Since the Kobayashi distance can only
decrease under compositions of holomorphic maps, we conclude that
$$k_{C^{\nu}}(\xi) \geq \epsilon_0 ||(i \circ \nu)_*
(\xi)||_{\omega}$$
for any tangent vector $\xi \in T_{C^{\nu}}$.
Integrating both sides of the inequality yields the proposition. $\Box$
\smallskip

Lang has conjectured that the converse also holds.  Let $X$ be a
compact, complex manifold endowed with a Hermitian metric.  Lang
conjectures that if there exists a positive constant $\epsilon$ such
that for every reduced, irreducible curve $C$ in $X$ the ratio of the
genus of the normalization of $C$ to the degree of $C$ (with respect
to the Hermitian metric) is bounded below by $\epsilon$, then $X$ is
hyperbolic. \smallskip

If the geometric genus of every curve in a variety $X$ is bounded
below by some fixed positive multiple of their degree, then $X$ does not
admit any non-constant holomorphic maps from any abelian variety. More
generally, Lang has conjectured that a projective variety $X$ is
hyperbolic if and only if it does not admit any holomorphic maps from
an abelian variety. The latter conjecture, of course, implies the
former one for projective varieties. \smallskip

One can ask for the relation between varieties of general type and
hyperbolic varieties. Since varieties of general type can
contain rationally connected subvarieties or abelian subvarieties, an
arbitrary variety of general type cannot be hyperbolic. However, Lang
has conjectured that the existence of subvarieties which are not of
general type accounts for the failure of hyperbolicity.

\begin{conjecture}\label{generaltype}
  A projective algebraic variety $X$ is hyperbolic if and only if
  every subvariety of $X$ is of general type.
\end{conjecture}   
Ein (\cite{ein:generaltype}) and later Voisin (\cite{voisin:clemens})
have shown that any subvariety of a very general hypersurface of
degree at least $2n+1$ in $\P^n$ is of general type. Combining their
theorem with this conjecture one obtains a conjectural sharp form of
Siu's Theorem \ref{higherdim}. \smallskip

Proving that a projective variety $X$ is hyperbolic usually has two
components.  One has to show that there are no non-algebraic maps of
$\C$ into $X$ and that there are no rational or elliptic curves in
$X$. When the problem is broken into these two components, then some
progress can be made on each component of the problem under some
geometric restrictions. We now survey some of the results on these
questions. \smallskip

One of the first people to make important progress on hyperbolicity
questions was the French mathematician Andr\'e Bloch, even though at
the time hyperbolicity was not defined
(see \cite{bloch:hyperbolicity}). Bloch begins by determining the
Zariski closure of an entire map into a complex torus.

\begin{theorem}\label{abelian}
  The Zariski closure of a holomorphic map $\C \rightarrow T$ to a
  complex torus $T$ is the translate of a subtorus of $T$.
\end{theorem}    
This theorem leads to a fairly complete understanding of the
hyperbolic subvarieties of complex tori.

\begin{corollary}
  A subvariety $X$ of an abelian variety $A$ which does not contain
  any translates of subtori of $A$ is hyperbolic.
\end{corollary}
Originally Bloch used these ideas to prove the following theorem often
referred to as Bloch's Theorem (\cite{bloch:hyperbolicity}).

\begin{theorem}[Bloch's Theorem]
  Any holomorphic map of $\C$ into a smooth, compact K\"ahler variety $X$
  whose irregularity ($h^0(X, \Omega_X^1)$) is bigger than its
  dimension is analytically degenerate, i.e. its image lies in a
  proper analytic subvariety.
\end{theorem}
{\bf Proof:} Recall that given a smooth, compact K\"ahler variety $X$, we
can associate to it a complex torus $\mbox{Alb}(X)$ of dimension
$h^0(X, \Omega_X^1)$ and a map $a: X \rightarrow \mbox{Alb}(X)$ with
the following universal property: for any complex torus $T$ and any
morphism $f: X \rightarrow T$, there exists a unique morphism $g:
\mbox{Alb}(X) \rightarrow T$ such that $f= g \circ a$. The complex
torus $\mbox{Alb}(X)$ is referred as the {\it Albanese}  variety of
$X$ and the
map $a$ is called the {\it Albanese map}.  $\mbox{Alb}(X)$ is  
unique up to isomorphism. \smallskip

Bloch's theorem follows by considering the Albanese map. Let $f: \C
\rightarrow X$ be a holomorphic map. Consider $a \circ f: \C
\rightarrow \mbox{Alb}(X)$.  Since the irregularity of $X$ is larger
than the dimension of $X$, the image of $X$ under the Albanese map
lies in a proper subvariety of $\mbox{Alb}(X)$.  By the universality
of the Albanese variety and its uniqueness, the Albanese image of $X$
cannot be the translate of a subtorus. Hence the image of $\C$ under
the map $a \circ f$ in the Albanese variety has to be analytically
degenerate in the image of $X$.  It follows that the image of the
original map $f: \C \rightarrow X$ has to be analytically degenerate.
$\Box$ \smallskip

Using these ideas one also obtains information about the hyperbolicity
of complements of ample divisors in abelian varieties.
\begin{theorem}\label{abeliancomplement}
  Let $A$ be an abelian variety. Let $D$ be an ample divisor that
  does not contain any translates of  abelian subvarieties.  Then
  $A-D$ is hyperbolic.
\end{theorem}

Recall that a holomorphic map from $\C$ into an algebraic variety is
called algebraically degenerate if the Zariski closure of the image
lies in a proper algebraic subvariety.  More recently McQuillan
(\cite{mcquillan:foliation}) has proved the algebraic degeneracy of
entire maps into surfaces of general type with $c_1^2 > c_2$.  More
precisely,

\begin{theorem}[McQuillan]
  If $X$ is a surface of general type satisfying $c_1^2 > c_2$, then
  all entire curves on $X$ are algebraically degenerate. In
  particular, if $X$ does not contain any rational or elliptic curves,
  then $X$ is hyperbolic.
\end{theorem}
McQuillan obtains this result by first proving a result about the
algebraic degeneracy of leaves of certain foliations on surfaces of
general type, then showing that under the assumptions on $X$ there
always is a foliation which contains the image of any entire map in
one of its leaves. \smallskip

Most ways of proving the algebraic degeneracy of entire maps into
complex projective manifolds depend on producing enough differential
relations on the manifold that every entire map has to satisfy. One
then shows that these relations have a small base locus. This forces
the holomorphic map to lie in this base locus. One often formalizes
these ideas in terms of jet bundles (see \cite{demailly:hyphyp},
\cite{demailly:hyperbolic} or any of Siu's papers cited above).
\smallskip

In the other direction, there has been extensive work on showing that
certain varieties do not contain any rational or elliptic curves. Here
we mention the work of G. Xu bounding below the geometric genus of any curve
on a very general surface of degree at least 5 in $\P^3$ (\cite{xu:bound},
see also Clemens' paper \cite{clemens:bound}).

\begin{theorem}[Xu]
On a very general surface of degree $d$ in $\P^3$, the geometric genus of
any curve is greater than or equal to $d(d-3)/2 - 2$. Tritangent plane
sections achieve this bound. For $d \geq 6$ the tritangent plane
sections are the only curves that achieve the bound.
\end{theorem}  
A very simple consequence of the theorem is that a general
hypersurface of degree at least 5 in $\P^3$ contains no rational or
elliptic curves. Thus the problem of showing the hyperbolicity of a
very general hypersurface in $\P^3$ of degree at least 5 reduces to
showing that any entire map into such a surface is algebraically
degenerate. \smallskip

Motivated by our discussion one can ask the following question:
\smallskip

\noindent {\bf Question:} Can the closure of the image of a
holomorphic map $\C \rightarrow \P^n$ be a variety of general type?
\smallskip

\noindent A negative answer to this question would imply one direction
of Conjecture \ref{generaltype}. If all the subvarieties of a compact
variety are of general type, then the variety would have to be
hyperbolic. At present this question seems very hard to answer.
\smallskip

We close this section with a discussion of how hyperbolicity varies in
families. Let $X$ be a compact $C^{\infty}$ manifold with Hermitian
metric $| \ |$. Brody proved that if one considers the various complex
structures that one can put on $X$, the set of those that are
hyperbolic is open in the analytic topology. More precisely, one can
consider the function
$$D(s) = \sup_{f \in \mbox{Hol} (\Delta, X_s)} |f'(0)|$$ on the moduli
spaces of complex structures parameterized by $S$. Brody proves

\begin{theorem}
  $D(s)$ is a continuous function. Since the hyperbolic complex
  structures correspond to those for which $D(s)$ is finite, the
  hyperbolic complex structures form an open set in the analytic
  topology.
\end{theorem}
Note that if one considers an algebraic family of quasiprojective
smooth varieties, then Brody's theorem may fail. In fact in
such a family (at least if we do not fix the $C^{\infty}$ type) the
set of fibers 
that are Brody  hyperbolic can be Zariski locally closed. Varying two of the
lines that meet a third at the same point in Figure \ref{lines}
provides such an example. \smallskip

Brody's theorem naturally raises the following question. \smallskip

\noindent {\bf Question:} Is hyperbolicity Zariski open in algebraic
families of projective varieties? \smallskip

\noindent  This question, to the best of my knowledge, is still
open. A positive answer would have interesting geometric implications.
For example, it is not known whether the space of high degree surfaces
in $\P^3$ that contain rational and elliptic curves form an algebraic
variety or a countable union of subvarieties of the space of surfaces.
A positive answer to the Zariski openness of hyperbolicity would
settle this and similar questions.

\section{The arithmetic of Kobayashi hyperbolicity}

In this section we summarize some conjectures, mainly due to Lang,
about rational points on hyperbolic varieties. We also give some
surprising consequences of the conjectures about the distribution of
rational points on curves. The main references for this section are
\cite{lang:hyperbolic}, \cite{caporaso:lang},
\cite{capjoemazur:points} and \cite{capjoemazur:lang}.

The arithmetic of curves of genus two or more and those of genus one
and zero exhibit drastically different properties. If $C$ is a curve
of genus zero or one defined over a number field $K$, then after
passing to a finite field extension $L$ there are infinitely many
$L$-rational points on $C$. In fact, $L$ can be chosen so that these
points are dense in the analytic topology. The Mordell-Faltings
Theorem stands in sharp contrast to this result.

\begin{theorem}[Mordell-Faltings]
  Let $C$ be a smooth curve of genus $g \geq 2$ defined over a number
  field $K$. Then $C$ has only finitely many rational points over any
  finite field extension $L$ of $K$.
\end{theorem}

\begin{definition} We say a variety $V$ defined over a number field
  $K$ is of {\it Mordell type} if $V$ contains only finitely many
  rational points over any finite field extension of $K$.
\end{definition}
  
It is natural to ask which, if any, of the geometric properties we
described so far imply that a variety is of Mordell type.  Clearly a
variety of Mordell type does not contain any rational or elliptic
curves. Similarly any map from an abelian variety to a variety of
Mordell type has to be constant. In view of these observations it is
not unreasonable to formulate the following conjecture due to Lang.

\begin{conjecture}\label{mordellic}
  A complex projective variety $V$ defined over a number field is of
  Mordell type if and only if it is hyperbolic.
\end{conjecture}

In fact Lang has conjectured much more precise statements. 

\begin{conjecture}[weak form]\label{weak}
  If $X$ is a variety of general type defined over a number field $K$,
  then the set of $K$-rational points of $X$ is not Zariski dense.
\end{conjecture}
The conjecture that a hyperbolic manifold is of Mordell type follows
from Conjecture \ref{weak} and the geometric Conjecture
\ref{generaltype} by the following argument.  If a hyperbolic manifold
$X$ has infinitely many rational points, then the Zariski closure of
these points is a subvariety not of general type by Conjecture
\ref{weak}. By Conjecture \ref{generaltype} every subvariety of $X$ is
of general type. Lang has conjectured an even stronger statement.

\begin{conjecture}[strong form]\label{strong}
  If $X$ is of general type, then there exists a proper algebraic
  subset $Z$ of $X$, such that over any finite extension $L$ of $K$,
  the number of $L$-rational points of $X-Z$ is finite.
\end{conjecture}

The Lang conjectures are currently open except for the case of curves
and more generally for subvarieties of abelian varieties. If true,
they would provide a fundamental understanding of the arithmetic of
varieties of general type.  \smallskip

The Lang Conjectures \ref{weak} and \ref{strong} have some surprising
consequences for rational points on curves of genus at least two.
These consequences were investigated by Caporaso, Harris and Mazur in
the papers cited above.

\begin{theorem}
  The weak form of Lang's conjecture implies that for every number
  field $K$ and genus $g \geq 2$, there exists a constant $B(K,g)$
  such that any curve of genus $g$ defined over $K$ has at most
  $B(K,g)$ $K$-rational points.
\end{theorem}

The strong form of Lang's conjecture implies an even more surprising
bound.
\begin{theorem}
  The strong Lang conjecture implies that for any $g \geq 2$ there
  exists an integer $N(g)$ such that there are only finitely many
  curves defined over a number field $K$ that have more than $N(g)$
  $K$-rational points.
\end{theorem}

\noindent {\bf Ideas behind the proof:} These theorems depend on the
  following geometric theorem Caporaso, Harris and Mazur prove.

\begin{theorem}(Correlation) \label{correlation}
  Let $f:X \rightarrow B$ be a proper morphism of irreducible and
  reduced schemes whose general fiber is a smooth curve of genus at
  least two. Then for $n$ sufficiently large the fiber product $h:
  X_B^n \dashrightarrow W$ admits a dominant rational map to a variety
  of general type. Moreover, if $f:X \rightarrow B$ is defined over a
  number field $K$, $h$ and $W$ can also be defined over $K$.
\end{theorem}
D. Abramovich proved a generalization of the Correlation Theorem for
morphisms with higher dimensional fibers of general type under
suitable hypotheses (see \cite{abramovich:correlation}). However, in
order to deduce the bounds on the rational points on curves, we only
need the case of curves. \smallskip

Caporaso, Harris and Mazur deduce the bounds on rational points on
curves by applying the Lang conjectures to a global family. More
precisely, they start with a family of curves $f : X \rightarrow B$
such that for every curve $C$ defined over $K$, there exists a
$K$-rational point on $B$ such that the fiber over it is isomorphic to
$C$ over $K$. \smallskip

Using Theorem \ref{correlation}, $X_B^n$, the $n-$th fiber product of
$X$ over $B$, admits a rational map to a variety $W$ of general type.
By Lang's conjecture there exists a smallest closed, proper subvariety
$V$ of $X_B^n$ that contains all the $K$-rational points. \smallskip

By studying the successive images of the complement of $V$ under the
projections to various factors of $X_B^n$, they produce a non-empty
open subset $U$ of $B$ and an integer $N$ such that for every rational
point $b \in U$ the fiber over $b$ has at most $N$ points. Then by
Noetherian induction on the complement of $U$, they conclude the
uniform bound. $\Box$

\bibliographystyle{math}
\bibliography{math}
\end{document}